\input amstex
\documentstyle{amams} % input Annals of Mathematics macros.
\document
\annalsline{156}{2002}
\received{October 23, 2000}
\startingpage{867}
\font\tenrm=cmr10
\def\joinrel{\mathrel{\mkern-4mu}}
\def\relbar{\mathrel{\smash-}}
\def\lrar{\relbar\joinrel\relbar\joinrel\relbar\joinrel\rightarrow}
\def\llrar{\relbar\joinrel\relbar\joinrel\relbar\joinrel\relbar\joinrel\relbar\joinrel\relbar\joinrel\relbar\joinrel\rightarrow}
\def\lmrar{\relbar\joinrel\relbar\joinrel\relbar\joinrel\relbar\joinrel\relbar\joinrel\rightarrow}
\def\rlrar{\relbar\joinrel\relbar\joinrel\relbar\joinrel\relbar\joinrel\relbar\joinrel\relbar\joinrel\relbar\joinrel\relbar\joinrel\relbar\joinrel\relbar\joinrel\rightarrow}

\def\mot#1{\,{\displaystyle\mathop{\to}^{#1}}\,}
\def\mol#1{\,{\displaystyle\mathop{\longrightarrow}^{#1}}\,}
\def\moll#1{\,{\displaystyle\mathop{\llrar}^{#1}}\,}
\def\morl#1{\,{\displaystyle\mathop{\rlrar}^{#1}}\,}

\def\molm#1{\,{\displaystyle\mathop{\lmrar}^{#1}}\,}
\catcode`\@=11
\font\twelvemsb=msbm10 scaled 1100

%\font\ninemsb=msbm7 scaled 1100%msbm9
\font\ninemsb=msbm10 scaled 800
\newfam\msbfam
\textfont\msbfam=\twelvemsb  \scriptfont\msbfam=\ninemsb
  \scriptscriptfont\msbfam=\ninemsb
\def\msb@{\hexnumber@\msbfam}
\def\Bbb{\relax\ifmmode\let\next\Bbb@\else
 \def\next{\errmessage{Use \string\Bbb\space only in math
mode}}\fi\next}
\def\Bbb@#1{{\Bbb@@{#1}}}
\def\Bbb@@#1{\fam\msbfam#1}
\catcode`\@=12

 \catcode`\@=11
\font\twelveeuf=eufm10 scaled 1100
\font\teneuf=eufm10
\font\nineeuf=eufm7 scaled 1100%eufm9
\newfam\euffam
\textfont\euffam=\twelveeuf  \scriptfont\euffam=\teneuf
  \scriptscriptfont\euffam=\nineeuf
\def\euf@{\hexnumber@\euffam}
\def\frak{\relax\ifmmode\let\next\frak@\else
 \def\next{\errmessage{Use \string\frak\space only in math
mode}}\fi\next}
\def\frak@#1{{\frak@@{#1}}}
\def\frak@@#1{\fam\euffam#1}
\catcode`\@=12

%--------------- Author macros ---------------

\define\br{\Bbb R}

\define\ci{\Cal I}
\define\cj{\Cal J}

\define\cl{\operatorname{cl}}
\define\cyl{\operatorname{cyl}}
\define\ocyl{\overset\circ\to\cyl}
\define\diam{\operatorname{diam}}
\define\hos{\operatorname{{holink}_s}}
\define\ho{\operatorname{holink}}
\define\id{\operatorname{id}}

\define\inclusion{\operatorname{inclusion}}

\def\ho{\hbox{\rm holink}}

\define\sing{\text{\rm sing}}

%-------------- Author entries --------------------
\title{The approximate tubular\\ neighborhood theorem} %Article title
\shorttitle{The approximate tubular  neighborhood theorem} % Shortened version for headline title
 \acknowledgements{Supported in part by NSF Grant DMS-9971367.\hfill\break
\hglue 28pt 1991 {\it Mathematics Subject Classification}.  Primary 57N80, 57N40; Secondary 55R65,
58A35.\hfill\break
\hglue 28pt {\it Key words and phrases}.  manifold stratified space, approximate tubular neighborhood,
tear\-drop neighborhood, manifold stratified approximate fibration.}
 \author{Bruce Hughes}
   \institutions{Vanderbilt University, Nashville,
Tennessee\\
 {\eightpoint {\it E-mail address\/}: bruce.hughes@vanderbilt.edu}}

\centerline{\tenpoint{\it In memory of James Robert Boyd Jr}\. (1921--2001)}

\vglue12pt 
\centerline{\bf Abstract}
\vglue12pt
Skeleta and other pure subsets of manifold stratified spaces are shown to have
neighborhoods which are teardrops of stratified approximate fibrations
(under dimension and compactness assumptions).
In general, the stratified approximate fibrations cannot be
replaced by bundles, and the teardrops cannot be replaced by
mapping cylinder neighborhoods. Thus, this is the best possible 
topological tubular neighborhood theorem in the stratified setting.

\vglue24pt

\section{Introduction}

One of the most striking differences between smooth and
topological manifolds concerns neighborhoods of submanifolds.  
For smooth manifolds
there is the classical Tubular Neighborhood Theorem of Whitney
asserting that every smooth submanifold has a neighborhood which
is the mapping cylinder of a smooth spherical fibre bundle. For
locally flat topological submanifolds, the examples of Rourke and
Sanderson \cite{26} show that neighborhoods which are mapping
cylinders of topological spherical fibre bundles need not exist.
However, Edwards \cite{7} proved that locally flat topological
submanifolds of manifolds of dimension greater than five do have
mapping cylinder neighborhoods, but the maps are a weak type of
bundle now called a manifold approximate fibration (see
\cite{20}).

For stratified spaces, there is a similar, but even more
pronounced, difference between the smooth and topological
categories. On the one hand, there are the smoothly stratified
spaces originally studied by Mather, Thom and Whitney (see
\cite{22}, \cite{29}, \cite{9}). Skeleta have mapping cylinder
neighborhoods whose maps are systems of topological fibre bundles
(see \cite{23}). On the
other hand, there are the topologically stratified spaces of
Siebenmann \cite{27} and Quinn \cite{25}. In this
setting skeleta (and even strata) may fail to have mapping
cylinder neighborhoods, and even when they do (as is the case for
locally flat submanifolds), the maps need not be fibre bundles.

The main result of this paper provides a substitute for the
missing mapping cylinder neighborhoods in topologically
stratified spaces. 

We work with the manifold
stratified spaces of Quinn \cite{25}.
These spaces are more general
than the locally conelike spaces of Siebenmann \cite{27}
in that Quinn's spaces are only locally conelike
up to stratified homotopy equivalence. In fact, the beauty of
Quinn's spaces is that their defining conditions are
homotopy-theoretical (from which geometric-topological properties
can be deduced). One point compactifications of manifolds with a
finite number of tame ends are examples of Quinn stratified
spaces which are locally conelike if and only if the manifolds
admit boundary completions. For another illustration of the
ubiquity of stratified spaces in the sense of Quinn, 
Cappell and Shaneson \cite{1} have
shown that mapping cylinders of stratified maps between smoothly
stratified spaces are manifold stratified spaces, 
even though they need not be smoothly stratified (cf\. \cite{13}).
For a survey on the various types of stratifications, as well as
related information, see Hughes and Weinberger \cite{21}.

Here is the main result of this paper.

\nonumproclaim{Approximate Tubular Neighborhood Theorem 1.1} 
Let $X$ be a manifold stratified space with compact singular set $X_{\sing}$
such that all the nonminimal strata of $X$ have dimension
greater than or equal to five. If $Y\subseteq X_{\sing}$ is
a pure subset of $X${\rm ,} 
then $Y$ has an approximate tubular neighborhood in $X$. 
\endproclaim 

The terminology in the theorem is explained fully in the sections
to follow, but here is a brief introduction. Pure subsets are
closed unions of strata, an important example being skeleta.
Approximate tubular neighborhoods are generalizations of
mapping cylinder neighborhoods of fibre bundles. Both
the mapping cylinder structure and the fibre bundle structure are
weakened. The mapping cylinder of a map $p:E\to B$ is replaced
by the teardrop of a map $q:U\to B\times\br$. A neighborhood of $B$
is constructed from this data by gluing $B$ to $U$ using the map
$q$. If $q$ were a fibre bundle, then this neighborhood would be
an open mapping cylinder of the desuspension of $q$ (in which the
$\br$ factor is split off). In general, $q$ is just required to
have a very weak homotopy lifting property, namely, $q$ is a
manifold stratified approximate fibration. Even though
desuspension is unavailable for these maps, there is still quite
a lot of geometry behind them.

The proof of Theorem 1.1 relies on both the statements and
techniques of special cases which have already been worked out.
First, there is the very important case of manifold stratified
spaces with only two strata studied by Hughes, Taylor, Weinberger
and Williams \cite{18}. Hughes and Ranicki \cite{17} specialized
further by requiring the lower stratum to be a point.
That single strata have approximate tubular neighborhoods was
established in \cite{15}.

The converse of Theorem 1.1, namely, that the teardrop of a
manifold stratified approximate fibration is a manifold stratified
space, was proved in~\cite{14}.
\pagegoal=49pc

Quinn indicated in his address to the International Congress
\cite{24} that topology with control is critical for the study of
singular and stratified spaces. Indeed, the
basic tools used in this paper come from controlled topology. In
particular, the geometric techniques have evolved from Chapman's
controlled engulfing methods \cite{3}, \cite{4}.
Other stratified tools come from \cite{12} and \cite{13}.
The main external input needed for the proof of Theorem 1.1 is
Quinn's Isotopy Extension Theorem \cite{25}.

Some applications of the approximate tubular neighborhood theorem
have already been outlined in the literature. Perhaps the most
important is the alternative approach to Weinberger's
surgery-theoretic classification of manifold stratified spaces
offered by him in  [30, p\. 189]  and  [32, pp\. 518--519].
The alternative approach applies to unstable classification
directly whereas Weinberger's first proof involves
stabilization-destabilization. 

Weinberger, in his address to the International Congress
\cite{31}, mentions applications to equivariant versions of
local contractiblity of homeomorphism groups and cell-like
approximation theorems. These results were first established by
Siebenmann \cite{27} and Steinberger and West \cite{28},
respectively, in the locally linear case. Another important
application is to complete the realization part of Quinn's
$h$-cobordism theorem \cite{25}. This was done for two strata by
Hughes, Taylor, Weinberger and Williams in \cite{18}. Further
applications, including multiparameter isotopy extension theorems
and Thom's isotopy lemmas, are mentioned in \cite{11}. Complete
details of these applications, along with results concerning
uniqueness, will be forthcoming.

This paper is organized as follows. Sections 2, 3 and 4 contain
the basic definitions and background information on manifold
stratified spaces, stratified approximate fibrations,
teardrops and approximate tubular neighborhoods.
Section 5 contains a special case which will
be used in the proof of the main result: it is shown that
collars of strata have approximate tubular neighborhoods. Section
6 establishes approximate tubular neighborhoods for certain
subsets of the singular set, the singular up-sets. Finally,
Section 7 contains the proof of the main result.

I want to thank Bruce Williams for the many stimulating
discussions we have had about the results in this paper.

%%%%%%%%%%%%%%%%%%%%%%%%%%%%%%%%%%%%%%%%%%%%%%%%%%%%%%%%%%%%%%%%%%%%%%%%%
%%%%%%%%%%%%%%%%%%%%%%%%%%%%%%%%%%%%%%%%%%%%%%%%%%%%%%%%%%%%%%%%%%%%%%%%%%%%%%
\section{Manifold stratified spaces}

This section contains the basic definitions from the theory of stratifications
as presented in \cite{11}, \cite{12}, \cite{13}, \cite{14},
\cite{18}, \cite{25}.

\pagegoal=48pc

\numbereddemo{Definition}
A {\it stratification\/} of a space $X$ consists of an index set $\ci$ and
a locally finite partition
$\{ X_i\}_{i\in\ci}$ of locally closed subspaces of~$X$ (the $X_i$ are pairwise disjoint and their union is $X$). For
$i\in\ci$,
$X_i$ is called the $i$-{\it stratum\/}
and the closed set
$$X^i = \mathbold{\cup}\{X_k ~|~ X_k\cap\cl (X_i)\not=\emptyset\}$$
is called the $i$-{\it skeleton\/}.
We say $X$ is a {\it space with a stratification\/}.
\enddemo

For a space $X$ with a stratification $\{ X_i\}_{i\in\ci}$, define
a relation $\leq$ on the index set $\ci$ by
$i\leq j$ if and only if $X_i\subseteq\cl (X_j).$
The {\it Frontier Condition\/} is satisfied 
if for every $i,j\in\ci$,
$X_i\cap\cl(X_j)\not=\emptyset$ implies $X_i\subseteq
\cl(X_j)$, in which case $\leq$ is
a partial ordering of $\ci$ and $X^i=\cl(X_i)$ for each $i\in\ci$. 

If $X$ is a space with a stratification satisfying the Frontier
Condition and $Y$ is a union of strata of $X$, then
$\cl(Y)\setminus Y$ is closed in $X$.

If $X$ is a space with a stratification, then a map $f:Z\times A\to X$ is
{\it stratum preserving along $A$\/} if for each $z\in Z$,
$f(\{ z\}\times A)$ lies in a single stratum of~$X$. In particular,
a map $f:Z\times I\to X$ is a {\it stratum preserving homotopy\/}
if $f$ is stratum preserving along $I$.
A homotopy $f:Z\times I\to X$ whose restriction to $Z\times [0,1)$ is
stratum preserving along $[0,1)$ is said to be {\it nearly stratum
preserving\/}.

\numbereddemo{Definition} \hskip-8pt
Let $X$ be a space with a stratification
$\{X_i\}_{i\in\ci}$ and $Y\subseteq X$.

 \vglue4pt
(1) $Y$ is {\it forward tame\/} in $X$
if there exist a neighborhood $U$ of $Y$ in $X$ and a 
homotopy $h:U \times I\to X$ such that $h_0 = \inclusion:U\to X$,
$h_t|Y = \inclusion:Y\to X$ for each $t\in I, h_1(U) = Y$, and 
$h((U\setminus Y)\times [0,1)) \subseteq X\setminus Y.$

 \vglue4pt
(2) The {\it homotopy link\/} of $Y$ in $X$ is defined by
$$\ho(X,Y)=\{ \omega \in X^I ~|~ \omega(t) \in Y\text{ if and only if }t=0\}.$$

\vglue4pt
(3) Let $x_0\in X_i\subseteq X$. The {\it local holink \/} (or
{\it local homotopy link}) {\it at \/}
$x_0$ is 
$$\multline
\ho(X,x_0) 
=\{\omega\in\ho(X,X_i) ~|~ \omega(0)=x_0  \\ \text{ and
} \omega(t)\in X_j \text{ for some } j, \text{ for all } t\in (0,1]\}.
\endmultline$$
\enddemo

All path spaces are given the compact-open topology.
Evaluation at $0$ defines a map
$q:\ho(X,Y)\to Y$  called {\it holink evaluation\/}.
There is a natural stratification of $\ho(X,x_0)$ into disjoint subspaces
$$\ho(X,x_0)_j = \{\omega\in\ho(X,x_0) ~|~ \omega(1)\in X_j\}.$$

\numbereddemo{Definition}
A space  $X$ with a stratification satisfying the Frontier Condition
is a {\it manifold  stratified space\/}
if the following four conditions are satisfied:

\vglue4pt
(1) {\it Forward tameness.} For each $k>i$, the stratum $X_i$ is
forward tame in $X_i \cup X_k$.

\vglue4pt
(2) {\it Normal fibrations.} For each $k>i$, the holink evaluation
$$q:\ho(X_i\cup X_k,X_i)\to X_i$$ 
is a fibration.

\vglue4pt
(3) {\it Compactly dominated local holinks.} For each
$x_0\in X$,
there exist a compact subset $C$
of the local homotopy link $\ho(X,x_0)$  and a stratum preserving
homotopy $$h:\ho(X,x_0)\times I\to\ho(X,x_0)$$ such that $h_0 =
{\rm id}$ and $h_1(\ho(X,x_0))\subseteq C$.

\vglue4pt
(4) {\it Manifold strata property.} 
$X$ is a locally compact, separable metric
space, each stratum $X_i$ is a topological manifold
(without boundary) and $X$ has only finitely many nonempty strata.
\enddemo

If $X$ is only required to satisfy conditions (1) and (2), then
$X$ is a {\it homotopically stratified space}.

\numbereddemo{Definition} \hglue-6pt
  The {\it singular set\/} of a space
$X$ with a stratification $\{X_i\}_{i\in\ci}$
is $$X_{\sing} = \mathbold{\cup}\{ X_i ~|~ \text{ for some } j\in\Cal I,
j\not= i,
\cl(X_j)\cap X_i \not=\emptyset\}.$$
In other words, $X_{\sing}$ is the union of all nonmaximal strata of $X$.
\enddemo

\numbereddemo{Definition} A subset $A$ of a space $X$ with a
stratification is a {\it pure subset\/} if $A$ is closed and is a
union of strata of $X$.
\enddemo

\vglue-6pt
\section{Stratified approximate fibrations}

We now give the definitions of the types of maps which are
important for manifold stratified spaces. 

\numbereddemo{Definition}
Let $X$ and $Y$ be spaces with stratifications $\{ X_i\}_{i\in\ci}$ and 
$\{ Y_j\}_{j\in\cj}$, respectively, and let $p:X\to Y$ be a map.

\vglue4pt
(1) $p$ is a {\it stratified fibration\/}
provided that given any space $Z$ and any commuting diagram
\pagebreak
$$\matrix
Z &{\displaystyle\mathop{\lrar}^f}& X \\
 {\scriptstyle\times 0}\Big\downarrow&&\Big\downarrow{\scriptstyle p} \\
Z\times I & {\displaystyle\mathop{\lrar}^{F}}& Y 
\endmatrix
$$
with $F$ a stratum preserving homotopy, there exists a
{\it stratified solution\/}; i\.e\., a stratum preserving homotopy
$\tilde F : Z\times I\to X$
such that $\tilde F (z,0) = f(z)$ for each
$z\in Z$ and $p\tilde F= F$. The diagram above is a {\it stratified
homotopy lifting problem\/}.

\vglue4pt
(2) $p$ is a {\it weak stratified approximate fibration\/}
provided that given any stratified homotopy lifting problem,
there exists a
{\it weak stratified controlled solution\/}; i\.e\., a map 
$\tilde F : Z\times I\times [0,1) \to X$
which is stratum preserving along  $I \times [0,1)$
such that $\tilde F (z,0,t) = f(z)$ for each
$(z,t) \in Z \times [0,1)$ and the function
$\bar F : Z \times I \times I\to Y$
defined by $\bar F|Z \times I \times [0,1) = p\tilde F$
and $\bar F|Z \times I \times \{ 1 \} = F \times\id_{\{1\}}$
is continuous.

\vglue4pt
(3) $p$ is a {\it manifold stratified approximate fibration\/} ({\rm MSAF}) if
$X$ and $Y$ are manifold stratified spaces and $p$ is a proper weak
stratified approximate fibration.

\vglue4pt
(4) If $\alpha$ is an open cover of $Y$, then $p$ is a {\it
stratified $\alpha$-fibration\/} provided that given any stratified
homotopy lifting problem,  there exists a
{\it stratified $\alpha$-solution\/}; i\.e\., a stratum preserving homotopy
$\tilde F : Z\times I\to X$
such that $\tilde F (z,0) = f(z)$ for each
$z\in Z$ and $p\tilde F$ is $\alpha$-close to $F$.

\vglue4pt
(5) $p$ is a {\it manifold approximate fibration\/} (MAF) if
$p$ is a {\rm MSAF} and $X$ and $Y$ have only one stratum each (i\.e\.,
they are manifolds).
\enddemo

See \cite{16} for clarification about weak stratified approximate
fibrations and how the definition above relates to definitions in
previous papers.

%%%%%%%%%%%%%%%%%%%%%%%%%%%%%%%%%%%%%
\vglue-6pt
\section{Teardrops and approximate tubular neighborhoods}

This section contains a review of the basic teardrop construction.
Given spaces $X$, $Y$ and a map $p:X\to Y\times\br$,
the {\it teardrop} of $p$ is the space denoted by 
$X\cup _p Y$ whose underlying set is the disjoint union
$X\amalg Y$ with the minimal topology such that
\vglue4pt
\item{(1)} $X\subset X\cup_p Y$ is an open embedding, and
\vglue4pt \item{(2)} the function $c: X\cup_p Y\to Y\times (-\infty ,
+\infty ]$ defined by
$$c(x)~=~\cases
p(x), &\text{if $x\in X$}\\
(x,+\infty ), &\text{if $x\in Y$~}
\endcases$$
is continuous.
\pagebreak

This is a generalization of the construction of the open
mapping cylinder of a map $g: X\to Y$. Namely, $\ocyl(g)$
is the teardrop $(X\times \br ) \cup_{g\times {\rm id}} Y$.
However, not all teardrops are open mapping cylinders because
not all maps to $Y\times\br$ can be split as a product. See
\cite{18} for more about the teardrop construction. 

If $X$ is a space with a stratification and $A\subseteq X$, 
we say $A$ has an {\it approximate tubular neighborhood in\/}
$X$ if
there is an open neighborhood $U$ of $A$ and an {\rm MSAF} 
$$p:U\setminus A\to A\times\br$$ 
such that the natural function $(U\setminus
A)\cup_p A\to U$ is a homeomorphism.
This has previously been called an {\rm MSAF} {\it teardrop neighborhood
in\/} X.
The condition is equivalent to saying
that $p$ is an {\rm MSAF} and the natural extension 
$$\tilde p:U\to A\times (-\infty,+\infty]$$ 
is continuous. In this case, $\tilde p$ is
also an {\rm MSAF} when $A\times (-\infty,+\infty]$ is given the
natural stratification (see [14, Prop\. 7.1], \cite{18}).

If $A$ does have an approximate tubular neighborhood in $X$, then
it is usually more convenient to replace $\br$ by $(0,+\infty)$
with $\{ 0\}\in [0,+\infty)$ playing the role of $\{ +\infty\}\in
(-\infty,+\infty]$. Thus, there is a map
of the form 
$$\varphi:U\to A\times [0,+\infty)$$ 
where $U$ is an open neighborhood of $A$ in $X$,
$\varphi^{-1}(A\times\{ 0\})= A$, $\varphi\vert:A\to A\times\{
0\}$ is the identity, and $\varphi$ is an {\rm MSAF}. 
This map 
$\varphi$ is called an 
{\it approximate tubular neighborhood map for\/} $A$ in $X$.

The following results show that the teardrop
construction yields manifold stratified spaces and that strata
have approximate tubular neighborhoods.

\proclaimtitle{\cite{14}}
\proclaim{Theorem}  If $X$ and $Y$ are manifold
stratified spaces each with only finitely many strata
and $p:X\to Y\times (0,+\infty)$ is a manifold stratified approximate
fibration{\rm ,}
then the teardrop $X\cup_pY$ with the natural stratification is a 
manifold stratified space.
\endproclaim

\proclaimtitle{\cite{15}}
\proclaim{Theorem} Let $X$ be a manifold stratified space with a
stratum $A$ satisfying\/{\rm :}
\vglue4pt
\item{\rm (1)} $A$ has compact closure $\cl(A)$ in $X${\rm ,}\vglue4pt
\item{\rm (2)} if $Z_1$ and $Z_2$ are distinct strata of $X$ 
with $Z_1\subseteq\cl(A)\cap\cl(Z_2)${\rm ,} then $\dim(Z_2)\geq 5$. 
\vglue4pt
\noindent Then $A$ has an approximate tubular neighborhood in $X$. 
\endproclaim
\pagebreak

%%%%%%%%%%%%%%%%%%%%%%%%%%%%%%%%%%%%%%%%%%%%%%%%%%%%%%%%%%%%%%%%%%%%%%%%%%
%%%%%%%%%%%%%%%%%%%%%%%%%%%%%%%%%%%%%%%%%%%%%%%%%%%%%%%%%%%%%%%%%%%%%%%%%%%%%%

\section{Approximate tubular neighborhoods of collars of strata}

In \cite{15} there is a proof that strata in manifold stratified
spaces have approximate tubular neighborhoods (under dimension
and compactness assumptions). In  this section
we extend that result slightly to show that collars of strata
have approximate tubular neighborhoods. This will be important in
the next section. 

The main result of this section is Proposition 5.2. Its proof
uses variations of the results in \cite{15} on ``stratified
sucking'' and ``homotopy near a stratum.'' The reader is required
to be familiar with those proofs. Lemma 5.1 shows how to deal
with the problem that limits of stratum preserving
processes need not be stratum preserving. 

Throughout this section let $N$ denote a manifold (without
boundary), possibly noncompact. Consider $N\times[0,+\infty)$ as
a manifold stratified space with two strata: $N\times\{ 0\}$ and
$N\times (0,+\infty)$. Let $\pi_1:N\times [0,+\infty)\to N$ and
$\pi_2:N\times [0,+\infty)\to [0,+\infty)$ denote the two projections.
We will assume that the one-point compactification of $N$ is a
manifold stratified space with two strata ($N$ and the point at
infinity) and that $\dim N \geq 5$ if $N$ is noncompact. This allows us to assume that $N$ has a metric with the
property that for every $\varepsilon >0$ there is a $\delta >0$ such
that any two maps into $N$ which are $\delta$-close are
$\varepsilon$-homotopic (rel any subset where the two maps agree).
Usually one would only have such a property for measurements made
by open covers. However, under the assumptions, $N$ has a
cocompact open subset which is the infinite cyclic cover of a
compact manifold \cite{17}. It follows that the desired metric
can be constructed. 
$N\times [0,+\infty)$ is given a standard product metric.

\proclaim{Lemma} Suppose $W$ is a manifold stratified space
and $p,p_n:W\to N\times [0,+\infty)$ are proper maps for
$n=1,2,3,\dots$ such that
\vglue4pt
\item{\rm (1)} $p_n$ is stratum preserving $1\over 2^n$\/{\rm -}\/homotopic to
$p_{n+1}$ for $n=1,2,3,\dots${\rm ,}
\vglue4pt
\item{\rm (2)} $p_n$ is a stratified $1\over 2^n$\/{\rm -}\/fibration for
$n=1,2,3,\dots${\rm ,}
\vglue4pt
\item{\rm (3)} $p=\lim_{n\to\infty}p_n$  {\rm (}\/uniformly\/{\rm ).}\/
\vglue4pt
Define $q:N\times [0,+\infty)\to N\times [0,+\infty)$ by
$$q(x,s)=\cases (x,0) & \text{if $0\leq s\leq 10$}\\ 
(x,s-10) & \text{if $ s\geq 10$}.\\
\endcases$$
Then $qp:W\to N\times [0,+\infty)$ is an  {\rm MSAF}. 
\endproclaim

\demo{Proof} 
We begin with a general construction which will be used in the
proof. This so called $\ast$-construction is used to convert 
homotopies into stratum
preserving homotopies. Let $F:Z\times I\to N\times [0,+\infty)$
be a homotopy such that $|\pi_2F(z,t)-\pi_2F(z,0)| < 1$ for each
$(z,t)\in Z\times I$. Define $F^\ast:Z\times I\to N\times
[0,+\infty)$ by $\pi_1F^\ast=\pi_1F$ and 
$$ 
\pi_2F^\ast(z,t)= 
{\cases
\pi_2F(z,0) & \text{if $\pi_2F(z,0)\leq 5$}\\
(\pi_2F(z,0)-5)\cdot\pi_2F(z,t)\\
\hskip.5in +(6-\pi_2F(z,0))\cdot\pi_2F(z,0) 
& \text{if $5\leq\pi_2F(z,0)\leq 6$}\\ 
\pi_2F(z,t)
& \text{if $6\leq\pi_2F(z,0)$.}\\
\endcases} $$
One can verify the following properties:
\vglue4pt
\item{(1)} $\pi_1F^\ast=\pi_1F$,
\vglue4pt
\item{(2)} $F_0^\ast=F_0$,
\vglue4pt
\item{(3)} $F^\ast$ is stratum preserving,
\vglue4pt
\item{(4)} $F^\ast(z,t)=F(z,t)$ if $\pi_2F(z,0)\geq 6$,
\vglue4pt
\item{(5)} $qF^\ast=qF$.
\vglue4pt

In the course of the proof we will use the observation $\pi_1qp=\pi_1p$.
Also, let $\hat q:N\times [10,+\infty)\to N\times [0,+\infty)$ be
the restriction of $q$ and note that $\hat q$ is a homeomorphism.
Now let 
$$\matrix
Z & {\displaystyle \mathop{\lrar}^{f}}& W \\
{\scriptstyle\times 0}\Big\downarrow&& \Big\downarrow {\scriptstyle qp}  \\
Z\times I &{\displaystyle \mathop{\lrar}^{F}} &N\times [0,+\infty)
\endmatrix$$
be a stratified homotopy lifting problem for which we are
required to  find a weak stratified controlled solution. Since
the problem is stratified, 
$$\pi_2F(z,t) > 0 \iff \pi_2qpf(z) > 0\iff\pi_2pf(z) >10.$$
Define $F':Z\times I\to N\times [0,+\infty)$ by
$$ F'(z,t) =\cases
(\pi_1F(z,t),\pi_2pf(z)) & \text{if $\pi_2pf(z)\leq 10$}\\
(\pi_1F(z,t),\pi_2\hat q^{-1}F(z,t))=\hat q^{-1}F(z,t)& \text{if
$\pi_2pf(z)\geq 10$.}\\
\endcases$$
One can verify that $F'$ is continuous, $qF'=F$ and 
$$\matrix
Z &{\displaystyle \mathop{\lrar}^{f}}& W \\
{\scriptstyle\times 0}\Big\downarrow&& \Big\downarrow {\scriptstyle p}  \\
Z\times I &{\displaystyle \mathop{\lrar}^{F'}} & N\times [0,+\infty)
\endmatrix$$
is a stratified homotopy lifting problem.

We will first show that for every $\varepsilon > 0$, the original
problem has a stratified $\varepsilon$-solution. 
For each $n=1,2,3\dots$ let $G^n:W\times [{n-1\over n},{n\over
n+1}]\to N\times [0,+\infty)$ be a stratum preserving ${1\over
2^n}$-homotopy from $p_n$ to $p_{n+1}$. The $G^n$ piece together
to define a map $W\times [0,1)\to N\times [0,+\infty)$ which is
stratum preserving along $[0,1)$. Moreover, $G:W\times I\to
N\times [0,+\infty)$ defined by $G\vert (W\times  [{n-1\over n},{n\over
n+1}]) = G^n$ and $G(w,1)=p(w)$ is a homotopy (continuous, but
not necessarily stratum preserving). For each $n=1,2,3,\dots$
consider the homotopy 
$H^n:Z\times [-{1\over n},1]\to N\times [0,+\infty)$ defined by
$$H^n(z,t)=\cases
G(f(z),t+1) & \text{if $-{1\over n}\leq t\leq 0$}\\
F'(z,t) & \text{if $0\leq t\leq 1$.}\\
\endcases$$ 
Note that $H^n(z,-{1\over n})=G(f(z),\frac{n-1}{n})=p_nf(z)$.
The $\ast$-construction yields a stratified homotopy lifting
problem
$$\matrix
Z &{\displaystyle\mathop{\lrar}^{f}}& W \\
{\scriptstyle\times \{-{1\over n}\}}\Big\downarrow&&\Big\downarrow   {\scriptstyle p_n}  \\
Z\times [-{1\over n},1]&{\displaystyle\mathop{\lrar}^{(H^n)^\ast}}& N\times [0,+\infty)
\endmatrix$$
which therefore has a stratified ${1\over 2^n}$-solution $\tilde
H^n:Z\times [-{1\over n},1]\to W$. 
Given $\varepsilon >0$, $n$ can be chosen large and 
$\tilde H^n$ can be reparametrized (by
covering the interval $[-{1\over n},0]$ quite rapidly) to get a
homotopy $\hat H^n:Z\times I\to W$ which is a stratified
$\varepsilon$-solution of the original problem; that is, $qp\hat
H^n$ is close to $F$ with the closeness depending on $n$. This
follows from the observations
$$p_n\tilde H^n\sim (H^n)^\ast \Rightarrow 
p\tilde H^n \sim (H^n)^\ast \Rightarrow 
qp\tilde H^n\sim q(H^n)^\ast=qH^n\sim qF'=F$$ where ``$\sim$''
denotes closeness which is small depending on $n$.

Now to get a stratified controlled solution from the existence of
stratified $\varepsilon$-solutions (for every $\varepsilon >0$) one
follows the proof of the corresponding unstratified result
\cite{19, Lemma 12.11} using the $\ast$-construction as needed.
 \enddemo

For more notation, let $X$ be a manifold stratified space
containing\break $N\times [0,+\infty)$ so that $N\times \{0\}$ and
$N\times (0,+\infty)$ are strata of $X$.

\proclaim{Proposition} Suppose 
\vglue4pt
\item{\rm (1)} $N\times [0,+\infty)$ has compact closure in $X${\rm ,}
\vglue4pt
\item{\rm (2)} if $Y$ and $Z$ are distinct strata of $X$ with
$Z\subseteq\cl(N\times [0,+\infty))\cap\cl(Y)${\rm ,} then $\dim(Y)\geq
5$. 
\vglue4pt\noindent
Then $N\times [0,+\infty)$ has an approximate tubular
neighborhood in $X$.
\endproclaim

\demo{Proof} 
From this point on we ask the reader to be familiar with the
proofs in \cite{17}, \cite{18} and, especially, \cite{15} of the
special cases of the approximate tubular neighborhood theorem.
Let $Z =\cl(N\times
[0,+\infty))\setminus (N\times [0,+\infty))$. 
In order for the techniques of [15, \S\S6,7] to apply we need
to assume that $Z$ is a single point. As in \cite{15} we can
reduce to this case by passing to the quotient space $X/Z$.
Since $N\times [0,+\infty)$ is stratified forward tame in $X$
\cite{12},
there is a
neighborhood $U$ of $N\times [0,+\infty)$ in $X$ and a nearly \pagebreak
stratum preserving deformation of $U$ to $N\times [0,+\infty)$ rel
$N\times [0,+\infty)$.
This deformation\break induces a map into the open mapping cylinder of the
holink evaluation\break $\hos(X,N\times [0,+\infty))\to N\times
[0,+\infty)$ (see [15, \S6]). The open mapping
cylinder has  natural $[0,+\infty)$-coordinates and there
results a proper map 
$f:U \to N\times[0,+\infty)\times [0,+\infty)$
such that
\vglue4pt
\item{(1)} $f^{-1}(N\times [0,+\infty)\times\{ 0\})= N\times [0,+\infty)$ and
$f\vert:f^{-1}(N\times [0,+\infty)\times\{ 0\})\to N\times
[0,+\infty)\times\{ 0\}$ is the identity, and
\vglue4pt
\item{(2)} the map $f$ by virtue of factoring through the mapping
cylinder of a stratified fibration (via a homotopy equivalence
with good control) has good enough lifting properties that for
every open cover $\beta$ of $N\times [0,+\infty)\times
(0,+\infty)$, $f$ is properly homotopic rel $N\times [0,+\infty)$
to a map $\tilde f:U\to N\times [0,+\infty)\times [0,+\infty)$
such that $\tilde f\vert:W\to N\times [0,+\infty)\times
(0,+\infty)$ is a proper stratified $\beta$-fibration where
$W=U\setminus (N\times [0,+\infty))$. 
\vglue4pt\noindent 
Given an open cover $\alpha$ of $N\times [0,+\infty)\times
(0,+\infty)$, if $\beta$ is fine enough, then the techniques of
\cite{15}, which are consequences of engulfing, show that $\tilde
f\vert$ is $\alpha$-close to a map $p:W\to N\times
[0,+\infty)\times (0,+\infty)$ with the property that there
exists a sequence $\{p_n\}_{n=1}^\infty$ of  proper maps
$p_n:W\to N\times [0,+\infty)\times (0,+\infty)$ such that
\vglue4pt
\item{(1)} $p_n$ is stratum preserving $1\over 2^n$-homotopic to
$p_{n+1}$ for $n=1,2,3,\dots$,
\vglue4pt
\item{(2)} $p_n$ is a stratified $1\over 2^n$-fibration for
$n=1,2,3,\dots$,
\vglue4pt
\item{(3)} $p=\lim_{n\to\infty}p_n$ (uniformly).
\vglue4pt\noindent 
If $\alpha$ is chosen correctly, then $p$ extends continuously to a
map $\tilde p:U\to N\times [0,+\infty)\times [0,+\infty)$ via the
identity $N\times [0,+\infty)\to N\times [0,+\infty)\times \{0\}$.
Define $q:N\times [0,+\infty)\times (0,+\infty)\to N\times
[0,+\infty)\times (0,+\infty)$ by
$$q(x,s,t)=\cases (x,0,t) & \text{if $0\leq s\leq 10$}\\ 
(x,s-10,t) & \text{if $ s\geq 10$}.\\
\endcases$$
We now apply Lemma 5.1 to the current situation by incorporating
the $(0,+\infty)$ factor into $N$. We conclude that
$qp:W\to N\times [0,+\infty)\times (0,+\infty)$ is an {\rm MSAF}. We
want to apply Lemma 6.11 below in order to conclude 
that there exists a stratum preserving homeomorphism
of $X$ onto the teardrop\break $(X\setminus (N\times
[0,+\infty)))\cup_{qp} (N\times [0,+\infty))$ which restricts to
the identity on\break $N\times [0,+\infty)$. In order to use Lemma 6.11
we need to observe that\break $q|:N\times [0,+\infty)\to N\times
[0,+\infty)$ extends to a stratum preserving map of $U$
to itself which is a homeomorphism on the complement of $N\times
[0,+\infty)$.
This is a special case of how Quinn's Isotopy Extension Theorem
\cite{25} is used in the proof of Lemma 6.7 below. Since $N\times
[0,+\infty)$ has an approximate tubular neighborhood in the
teardrop, the proof is complete.
\enddemo
\pagebreak

%%%%%%%%%%%%%%%%%%%%%%%%%%%%%%%%%%%%%%%%%%%%%%%%%%%%%%%%%%%%%%%%%%%%%%%%%%%%%%
\section{Approximate tubular neighborhoods for singular up-sets}

This section establishes the result that certain subsets of the
singular set have approximate tubular neighborhoods. This will be
crucial in the inductive proof of the main result.

\numbereddemo{Definition} If $X$ is a space with a
stratification $\{ X_i\}_{i\in\ci}$, then a subset $Y$ of
$X_{\sing}$ is a {\it singular up-set\/} of $X$ if $Y$ is a union
of strata of $X_{\sing}$ and if  whenever $X_j$
is a stratum of $X_{\sing}$ for which there exists a stratum
$Y_i$ of $Y$ with $Y_i\subseteq \cl(X_j)$, then $X_j\subseteq Y$.
\enddemo

Note that a singular up-set $Y$ in $X$ need not be closed in $X$.

\proclaim{Theorem} Let $X$ be a manifold stratified space
with compact singular set $X_{\sing}$ and let $Y\subseteq
X_{\sing}$ be a singular up\/{\rm -}\/set of $X$ satisfying\/{\rm :}
\vglue4pt \itemitem{} if $Z_1$ and $Z_2$ are distinct strata of $X$ with
$Z_1\subseteq \cl(Y)\cap\cl(Z_2)${\rm ,}\hfill\break then $\dim(Z_2)\geq 5$. 
\vglue4pt\noindent
Then $Y$ has an approximate tubular neighborhood in $X$.  
\endproclaim

\demo{Proof} The proof is by induction on the number 
$k$ of strata of $Y$. If $k=1$, then $Y$ is a stratum
of $X$ with compact closure and the Main Theorem of \cite{15}
implies that $Y$ has an approximate tubular neighborhood in $X$.

Assume $k>1$ and Theorem 6.2 holds in the case of fewer than $k$
strata. Write $Y=A\cup B$ where $A$ is a minimal stratum of $Y$
and $B=Y\setminus A$. Using \cite{15} again, $A$ has an
approximate tubular 
neighborhood in $X$. By the inductive hypothesis, $B$
has an approximate tubular neighborhood in $X$. 
Of course, $A$ also has an approximate tubular neighborhood in $Y$
(by restricting the approximate tubular neighborhood map for
$A$ in $X$ to $Y$).

\vglue12pt {\it Step} 1  (\/{\it  Notation for the approximate tubular neighborhoods}\/).
Let $U_B$ be an open neighborhood of $B$ in $X$ for which
there is an approximate tubular neighborhood  map 
$\varphi_B:U_B\to B\times [0,+\infty)$.
Let $V_A$ be an open neighborhood of $A$ in $Y$ for which
there is an approximate tubular neighborhood map\break
$\varphi_A:V_A\to A\times [0,+\infty)$.
We need to show how to modify $\varphi_A$ so that it has the
additional property:
\vglue4pt\itemitem{}
if $x\in \cl(A)\setminus A$ and $U$ is an open neighborhood of
$x$ in $X$,\hfill\break then there exists an open neighborhood $V$ of $x$ in
$X$ such\hfill\break  that $\varphi_A^{-1}(\{ a\}\times [0,11])\subseteq U$
whenever $a\in V\cap A$.
\vglue4pt\noindent 
To this end let $\rho:\cl(A)\to I$ be a map such that
$\rho^{-1}(0)=\cl(A)\setminus A$\break and $\diam\varphi_A^{-1}(\{
a\}\times [0,\rho(a)])$ goes to $0$ as $a\in A$ approaches
$\cl(A)\setminus A$.\break Define $\rho':A\times [0,+\infty)\to  A\times
[0,+\infty)$ by
$\rho'(a,s)= (a,{11s\over\rho(a)})$, \pagebreak and\break
$\varphi_A':V_A\to A\times [0,+\infty)$ by
$\varphi_A'=\rho'\circ\varphi_A$. Then $\varphi_A'$ is the
approximate tubular neighborhood map with the additional property
(but we retain the notation $\varphi_A$ for this map).
\vglue12pt

{\it Step} 2 (\/{\it  Modifying $X$ along $V_A$}\/).
Let $A'$ be the one-point compactification of $A\times
[0,+\infty)$ with the point at infinity denoted $\omega$.

\proclaim{Claim} $A'$ is a manifold stratified space with strata
$\{\omega\}${\rm ,} $A\times \{ 0\}${\rm ,} $A\times (0,10)${\rm ,} $A\times \{
10\}$ and $A\times (10,+\infty)$.
\endproclaim

\demo{Proof} Note that $\cl(A)$ is a compact manifold stratified
space with $\cl(A)\setminus A$ as a closed manifold stratified
subspace. Stratify $[0,+\infty]$ with strata $\{ 0\}$, $(0,10)$,
$\{ 10\}$, $(10,+\infty)$ and $\{ +\infty\}$. Give $\cl(A)\times
[0,+\infty]$ the product stratification (which makes it a
manifold stratified space [16, 4\.1]). Since
$Z=((\cl(A)\setminus A)\times [0,+\infty])\cup(\cl(A)\times\{
+\infty\})$ is a compact manifold stratified subspace and
$A'=(\cl(A)\times [0,+\infty])/Z$, it 
follows from \cite{15} that $A'$ is a manifold stratified space.
 \enddemo

Define $\varphi_A':\cl(Y)\to A'$ by $\varphi_A'\vert
V_A=\varphi_A:V_A\to A\times [0,+\infty)\subseteq A'$ and
$\varphi_A'(\cl(Y)\setminus V_A)=\omega$.

\proclaim{Claim} $\varphi_A':\cl(Y)\to A'$ is a {\rm MSAF}.
\endproclaim

\demo{Proof}
Let
$$\matrix
Z &{\displaystyle \mathop{\lrar}^{f}}& \cl(Y) \\
{\scriptstyle\times 0}\Big\downarrow&&  \Big\downarrow {\scriptstyle\varphi_A'} \\
Z\times I &{\displaystyle\mathop{\lrar}^{F}}& A'
\endmatrix$$
be a stratified homotopy lifting problem. Let
$Z_\omega=f^{-1}(\cl(Y)\setminus V_A)$. Then the problem above
restricts to a stratified lifting problem
$$\matrix
Z\setminus Z_\omega &{\displaystyle \mathop{\lrar}^{f|}}&  V_A \\
{\scriptstyle\times 0}\Big\downarrow&&  \Big\downarrow {\scriptstyle\varphi_A} \\
(Z\setminus Z_\omega)\times I &{\displaystyle\mathop{\lrar}^{F|}}& A\times [0,+\infty)
\endmatrix$$
which has a stratified controlled solution $\tilde F:(Z\setminus
Z_\omega)\times I\times [0,1)\to V_A$. 
It is not too hard to modify $\tilde F$ so that $\diam\tilde F(\{
z\}\times I\times [0,1))$ goes to zero as $z$ gets close to
$Z_\omega$. If that modification is made, then $\tilde F$ will
extend to a stratified controlled solution of the original
problem by setting $\tilde F(z,s,t) =f(z)$ if $z\in Z_\omega$.
 \enddemo

Define the attaching space $X'=X\cup_{\varphi_A'}A'$. It follows
from [16, 6\.2] that $X'$ is a manifold stratified space.
Let $q_X:X\amalg A'\to X'$ be the quotient map.
 \vglue12pt
{\it Step} 3 (\/{\it  A neighborhood of the collar $A\times [0,10)$ 
in $X'$}\/).
Use Proposition 5.2 to get an open 
neighborhood 
$W_A$ of $A\times [0,10)$ in
$X'$ and a proper map 
$$\xi_A:W_A\to A\times
[0,10)\times [0,+\infty)$$ 
such that
$\xi_A^{-1}(A\times [0,10)\times\{0\})= A\times [0,10)$ and
$\xi_A|:A\times [0,10) \to A\times
[0,10)\times\{0\}$ is the identity.
\enddemo

\vglue12pt
{\it Step} 4 (\/{\it  Using uniqueness}\/).

\proclaim{Lemma} Let $M$ and $N$ be manifolds without
boundary such that $\dim(M)\geq 5$. Suppose $p:M\times I\to
N\times (0,10)\times I$ is a $1$\/{\rm -}\/parameter family of manifold
approximate fibrations\/{\rm ;} that is{\rm ,} $p$ is fiber preserving
over $I$ and $p_t:M\to N\times (0,10)$ is a manifold approximate
fibration for each $t\in I$.
Then there exists a manifold approximate fibration $\hat p:M\to
N\times (0,10)$ such that $\hat p=p_0$  over $N\times (0,2)$ and
$\hat p=p_1$ over $N\times (8,10)$. 
\endproclaim

\demo{Proof} By the straightening principle \cite{10} (cf\. \cite{19})
there is an isotopy $H:M\times I\to M\times I$ with $H_0=\id_M$
such that $pH$ is as close to $p_0\times\id_I$ as desired. By the
estimated homotopy extension theorem \cite{5}, there is a map
$\tilde p:M\to N\times (0,10)$ such that $\tilde p= p_0$ over
$N\times (0,4)$, $\tilde p= p_1H_1$ over $N\times (6,10)$, 
and $\tilde p$ is close to
$p_0$. By the sucking principle \cite{10} (cf\. \cite{19}), we may
additionally assume that $\tilde p$ is a manifold approximate
fibration. By the isotopy extension theorem \cite{8} there is a
homeomorphism $h:M\to M$ such that $h=\inclusion$ on $\tilde
p^{-1}(N\times (0,3))$ and $h=H_1$ on $\tilde p^{-1}(N\times (7,10))$.
Finally, $\hat p =\tilde ph^{-1}$ is the desired manifold
approximate fibration. 
 \enddemo

Returning to the proof of Theorem 6.2, let
$$U_{AB}=\varphi_B^{-1}(\varphi_A^{-1}(A\times (0,10))\times
[0,+\infty))\subseteq U_B \subseteq X$$ and define the composition
$$\varphi_{AB}:U_{AB} \mol{\varphi_B}  \varphi_A^{-1}(A\times
(0,10))\times [0,+\infty) \moll{\varphi_A\times\id_{[0,+\infty)}} 
A\times (0,10)\times [0,+\infty).$$ Note that $\varphi_{AB}$ is
an {\rm MSAF} as follows. First, $\varphi_A\times\id_{[0,+\infty)}$ is
an {\rm MSAF} [16, 4\.3]. Then [16, 7\.4, 4\.5] also implies that
$\varphi_{AB}$ is an {\rm MSAF}.
Now let $U_{AB}'= q_X(U_{AB})\subseteq X'$.
There is an induced map 
$$\varphi_{AB}':U_{AB}'\to A\times
(0,10)\times [0,+\infty).$$
That is, $\varphi_{AB}'\circ  q_X=\varphi_{AB}:
U_{AB}\to A\times (0,10)\times [0,+\infty)$.
Note that $\varphi_{AB}'$ has the following three properties:
\vglue4pt
\item{(1)} $(\varphi_{AB}')^{-1}(A\times (0,10)\times\{ 0\})= A\times (0,10)$,
\vglue2pt
\item{(2)} $\varphi_{AB}'\vert:A\times (0,10)\to A\times \{ 0\}$ 
is the identity,
\vglue2pt
\item{(3)} $\varphi_{AB}'\vert:U_{AB}\setminus (A\times (0,10))\to A\times (0,10)\times (0,+\infty)$ is an {\rm
MSAF}.
\vglue4pt\noindent 
The first two properties are obvious. For the third, note that since\break
$q_X\vert: U_{AB}\setminus (A\times (0,10))\to U_{AB}'
\setminus (A\times (0,10))$ is a homeomorphism, we can express
$\varphi_{AB}'$ on $U_{AB}\setminus (A\times (0,10))$ as $\varphi_{AB}
\circ q_X^{-1}$, and $\varphi_{AB}$ is an {\rm MSAF}. It follows from
\cite{14, Prop. 7.1} that $\varphi_{AB}':U_{AB}\to A\times (0,10)\times
[0,+\infty)$ is an~{\rm MSAF}.

Thus, we have two maps
$\xi_A|:W_A\setminus (A\times\{0\})\to A\times (0,10)\times[0,+\infty)$ and
$\varphi_{AB}':U_{AB}'\to A\times (0,10)\times [0,+\infty)$
which are {\rm MSAF}s 
and give
approximate tubular neighborhoods of $A\times (0,10)$ in $X'$.
Moreover, over $A\times (0,10)\times (0,+\infty)$ these {\rm MSAF}s are
actually MAFs (because their inverse images miss $X_{\sing}$). It
follows from the uniqueness results of \cite{18} and Lemma 6.5
that there exists a neighborhood $W_A'$ of $A\times [0,10)$ in
$X'$ and a map
$$\xi_A':W_A'\to A\times [0,10)\times
[0,+\infty)$$ 
such that 
\vglue4pt
\item{(1)} $W_A'\subseteq W_A\cup U_{AB}'$,
\vglue2pt
\item{(2)} $\xi_A'=\xi_A$ over $A\times (0,2)\times [0,+\infty)$,
\vglue2pt
\item{(3)} $\xi_A'=\varphi_{AB}'$ over $A\times (8,10)\times
[0,+\infty)$,
\vglue4pt
\item{(4)} $(\xi_A')^{-1}(A\times [0,10)\times\{0\})=A\times [0,10)$
and $\xi_A'|:A\times [0,10)\to A\times [0,10)\times\{0\}$ is the
identity,
\vglue4pt
\item{(5)} $\xi_A'$ is an {\rm MSAF} over $A\times (0,10)\times [0,+\infty)$.
\enddemo
 
{\it Step} 5 (\/{\it  Shrinking and pushing}\/).
Let $\hat X$ be the quotient space obtained from $X$ with
the equivalence relation generated by setting $x\sim y$
if $x,y\in V_A$, $\varphi_A(x)=(z,s)$ and $\varphi_A(y)=(z,t)$
for some $z\in A$ and $0\leq s,t \leq 10$. 
Let $g:X\to \hat X$ be the quotient map.
We identify $A$
with its image under $g$ so that $g|:A\to
A$ is the identity. 
Because of the additional condition imposed on $\varphi_A$ in
Step 1, it follows that $g$ is a closed map (which is to say that
the induced decomposition of $X$ is upper semicontinuous; cf\.
Daverman [6, p\. 8]). Moreover, $\hat X$ is a locally compact,
separable metric space and $g$ is a proper map [6, pp\. 13--17].

\proclaim{Claim} There exists a homeomorphism $\tilde
g:X\to\hat X$ such that
\vglue4pt
\item{\rm (1)} $\tilde g\vert:A\to A$ is the identity{\rm ,}
\vglue2pt
\item{\rm (2)} $\tilde g(Y)\subseteq g(Y)${\rm ,} and{\rm ,} in fact{\rm ,} $\tilde
g(S)\subseteq g(S)$ for each stratum $S$ of $Y${\rm ,}
\vglue2pt
\item{\rm (3)} if $S$ is a stratum of $X$ missing $Y${\rm ,} then $\tilde
g(S)=g(S)$.
\endproclaim

The proof of Claim 6.6 is based on Bing's Shrinking Criterion (cf\.
\cite{2}, \cite{6}) and 
the following lemma:

\proclaim{Lemma} There exist stratum preserving 
shrinking homeomorphisms for $g$\/{\rm ;} that is{\rm ,} for each open cover $\Cal U$ 
of $X$ and $\Cal V$ of $\hat X${\rm ,} there exists a stratum preserving 
homeomorphism $H:X\to X$  {\rm (}\/that is{\rm ,} if $S$ is a stratum of $X${\rm ,}
then $H(S)=S${\rm )} such that $gH$ is $\Cal V$\/{\rm -}\/close to $g${\rm ,} each 
$Hg^{-1}(y)$ lies in some element of $\Cal U$
and $H\vert\cl(A)$ is the inclusion.
\endproclaim

\demo{Proof} 
An isotopy $h_t, t\in I$, of $A\times
(0,+\infty)$ affecting only the $(0,+\infty)$-coordinates
and moving $A\times (0,10]$ close to $A\times\{0\}$ can be
approximately lifted to a stratum preserving isotopy $\tilde h_t, t\in I$ of
$V_A\setminus A$ so that $\varphi_A\tilde h_t$ is as close as
needed to $h_t\varphi_A$. 
This comes from using the engulfing result [15, 4.3] together
with Chapman's stacking technique [2, Lemma 3.5].
Now Quinn's Isotopy Extension Theorem
\cite{25} implies that $\tilde h_t$ can be extended to a stratum
preserving isotopy $\hat h_t$ of all of $X$. This extension is
done one stratum at a time in such a way that the desired control
is retained. Then the
$\hat h_1$ provide the required shrinking homeomorphisms
for $g$. 
 \enddemo

\demo{Proof of Claim {\rm 6.6}}
The proof of Bing's Shrinking Criterion given in
\cite{2} provides a proper map $k:X\to X$ constructed as a limit
$k=\lim_{n\to\infty} H_1\circ H_2\circ\cdots\circ H_n:X\to X$
where the $H_i$'s are shrinking homeomorphisms given by Lemma 6.7 
so that $\tilde g=g\circ k^{-1}$ defines the desired homeomorphism.
 \enddemo

Let $Y'$ be the quotient space obtained from $Y$ with the
equivalence relation generated by setting $x\sim y$ if $x,y\in
V_A$ and $x,y\in\varphi_A^{-1}(z,s)$ for some $(z,s)\in A\times
[0,10]$. Note that $Y'$ contains a natural copy of $A\times
[0,10]$. In fact, $Y'$ is the attaching space
$$Y'= Y\cup_{\varphi_A\vert}(A\times [0,10])$$ 
where $\varphi_A\vert:\varphi_A^{-1}(A\times [0,10])\to A\times
[0,10]$.
Let 
$$q_Y:Y\to Y'$$ 
be the quotient map. Write $Y'=(A\times [0,10])\cup B'$ where
$(A\times [0,10])\cap B'=\emptyset$ and $(A\times [0,10])\cap\cl(B')
=A\times\{ 10\}$.

Let $Y''$ be the quotient space obtained from $Y'$ with the
equivalence relation generated by setting $(x,s)\sim (x,t)$ for
each $x\in A$ and $0\leq s,t\leq 10$. Then $Y''\subseteq \hat X$.
Let $q_{Y'}:Y'\to \hat X$ be the composition of the quotient map
$Y'\to Y''$ followed by the inclusion $Y''\to \hat X$. 

Let $\tilde\beta:X\to X$ be the map $\tilde\beta=\tilde
g^{-1}\circ g$,
let $\beta =\tilde\beta|:Y\to Y$, and let 
$$\pi:Y'=(A\times [0,10])\cup B'\to Y=A\cup B$$ 
be the map
$\pi=\tilde g^{-1}\circ q_{Y'}$.
We call $\pi$ {\it the push\/}; it is the key geometric move which
allows the meshing of the two approximate tubular neighborhoods.
Note that $\pi\vert:A\times [0,10]\to A\subseteq Y$ is the projection
and $\pi\vert:\cl(B')\to Y$ is a stratum preserving homeomorphism.
Essentially, $\pi$ is the collapse of an external collar.

\proclaim{Claim} $Y'$ is a manifold stratified space with
strata $A\times\{ 0\}${\rm ,} $A\times (0,10)${\rm ,} $A\times\{ 10\}$ and
$q_Y(S)\cap B'$ for each stratum $S$ of $B\subseteq Y$.
\endproclaim

\demo{Proof} First $\cl(B')$ is a manifold stratified space
because $\pi\vert:\cl(B')\to Y$ is a stratum preserving
homeomorphism. Then $Y'$ is a manifold stratified space by the
adjunction theorem of [16, 6\.2].
 \enddemo

\proclaim{Claim} $\pi:Y'\to Y$ is a {\rm MSAF}.
\endproclaim

\demo{Proof} Let
$$\matrix
Z &{\displaystyle \mathop{\lrar}^{f}}& Y' \\
{\scriptstyle\times 0}\Big\downarrow&&\Big\downarrow {\scriptstyle\pi} \\
Z\times I &{\displaystyle \mathop{\lrar}^{F}}& Y
\endarray$$
be a stratified homotopy lifting problem. Define $\tilde
F:Z\times I\to Y'$ by $\tilde F(z,t) =
(\pi\vert\cl(B'))^{-1}F(z,t)$ if $f(z)\in \cl(B')$. If $f(z)\in
A\times [0,10]$, define $\tilde F(z,t)$ by setting $\pi_1\tilde
F(z,t)=F(z,t)$ and $\pi_2\tilde F(z,t)=\pi_2f(z)$ where
$\pi_1:A\times [0,10]\to A$ and $\pi_2:A\times [0,10]\to [0,10]$
are the projections. Then $\tilde F$ is a stratified solution
showing that $\pi$ is actually a stratified fibration.
 \enddemo

\proclaim{Claim} $q_Y:Y\to Y'$ is a {\rm MSAF}.
\endproclaim

\demo{Proof} This follows from [16, 7\.1].
 \enddemo

 {\it Step} 6 (\/{\it  Recognizing a teardrop}\/).
The plan is to define a neighborhood $U$ of $Y$ in $X$ together
with a {\rm MSAF} $\varphi:U\setminus Y\to Y\times (0,+\infty)$. This
map will not extend via the identity $Y\to Y\times \{0\}$ so
that we will not be able to conclude immediately that this gives
an approximate tubular neighborhood of $Y$ in $X$. However, $\varphi$ will extend
via $\beta$ and we will then be able to draw the necessary
conclusions from the following lemma (with $U$ playing the role
of $X$ so that $Y$ is closed).

\proclaim{Lemma} Let $X$ and $Y=A\cup B$ be as above{\rm ,} but now
assume that $Y$ is closed in $X$.
Suppose $\tilde\beta:X\to X$ is a proper surjection such
that\/{\rm : }
\vglue4pt
\item{\rm (1)} $\tilde\beta^{-1}(Y)=Y$ and $\beta:Y\to Y$ denotes the
restriction of $\tilde\beta${\rm ,} 
\vglue4pt
\item{\rm (2)} $\beta^{-1}(A)=N$ is a closed neighborhood of $A$ in $Y${\rm ,}
\vglue4pt
\item{\rm (3)} $\beta|:A\to A$ is the identity{\rm ,}
\vglue4pt
\item{\rm (4)} $\tilde\beta|:X\setminus N\to X\setminus A$ is a
homeomorphism.
\vglue4pt\noindent 
Suppose further that $\varphi:X\to Y\times [0,+\infty)$ is a
proper map such that\break $\varphi^{-1}(Y\times\{0\})=Y$ and
$\varphi(x)=(\beta(x),0)$ for each $x\in Y$.
Then there is a homeomorphism $h:X\to (X\setminus
Y)\cup_{\varphi|(X\setminus Y)}Y$ which restricts to the identity
on~$Y$. Moreover{\rm ,} if $\tilde\beta$ is stratum preserving in the
sense that $\tilde\beta(S)=S$ for each stratum $S$ of $X\setminus
Y${\rm ,} then $h$ also has this property.
\endproclaim

\demo{Proof} Define $h:X\to  (X\setminus Y)\cup_{\varphi|(X\setminus
Y)}Y$ by 
$$h(x)=\cases
x& x\in Y,\\
\tilde\beta^{-1}(x)& x\in X\setminus Y.\\
\endcases$$
Clearly, $h$ is a bijection. The continuity criterion from
\cite{18} can be used to see that $h$ is continuous as follows.
First, one
needs to check that $h|:X\setminus Y\to 
(X\setminus Y)\cup_{\varphi|(X\setminus Y)}Y$ is an open
embedding. But this map is $\tilde\beta^{-1}$, so this is
obvious. Second, letting $c:
(X\setminus Y)\cup_{\varphi|(X\setminus Y)} 
\to Y\times [0,+\infty)$ be the teardrop collapse, one must check 
that $c\circ h:X\to Y\times [0,+\infty)$ is continuous.
This map is seen to be 
$$x\mapsto\cases 
(x,0)& x\in Y,\\
\varphi\tilde\beta^{-1}(x)& x\in X\setminus Y.\\ 
\endcases$$
Let $x_n\in X\setminus Y$ $n=1,2,3,\dots$ be a sequence with
$x_n\to x_0\in Y$ and show that $\varphi\tilde\beta^{-1}(x_n)\to
(x_0,0)$. If $x_0\in B$, then
$\tilde\beta^{-1}(x_n)\to\tilde\beta^{-1}(x_0)=\beta^{-1}(x_0)$.
Thus, $\varphi\tilde\beta^{-1}(x_n)\to
\varphi\beta^{-1}(x_0)=(x_0,0)$. If, on the other hand, $x_0\in
A$, it follows from the local compactness of $X$ and the
properness of $\tilde\beta$ that after passing to a subsequence
we may assume that $\tilde\beta^{-1}(x_n)\to x_0'$ for some
$x_0'\in X$. Then $x_n\to\tilde\beta(x_0')$ and so
$\tilde\beta(x_0')=x_0\in Y$. Thus, $x_0'\in Y$ and
$\varphi(x_0')=(\beta(x_0'),0)=(x_0,0)$. Finally,
$\varphi\tilde\beta^{-1}(x_n)\to\varphi(x_0')=(x_0,0)$ as
desired.  Hence, $h$ is continuous.

To see that $h^{-1}$ is continuous, first note that it is given
by 
$$x\mapsto\cases
x& x\in Y,\\
\tilde\beta(x)& x\in X\setminus Y.
\endcases$$
It suffices to consider a sequence $x_n\in X\setminus Y$,
$n=1,2,3,\dots$ such that $x_n\to x_0\in Y$ in the teardrop of
$\varphi$ and show that $\tilde\beta(x_n)\to x_0$ in $X$. We know
that $\varphi(x_n)\to (x_0,0)$ in $Y\times [0,+\infty)$ (because
the teardrop collapse is continuous). By the local compactness of
$Y$ and the properness of $\varphi$, we may assume after passing
to a subsequence that $x_n\to x_0'$ for some $x_0'\in X$. Thus,
$\varphi(x_n)\to\varphi(x_0')$. So $\varphi(x_0')=(x_0,0)\in
Y\times\{ 0\}$ and so $\beta(x_0')=x_0$. Also $\varphi(x_n)\to
(x_0,0)$. Note $x_0'\in Y$. Now $x_n\to x_0'$ in $X$ implies
$\tilde\beta(x_n)\to\tilde\beta(x_0')$ in $X$ which in turn
implies $\tilde\beta(x_n)\to\beta(x_0')=x_0$ as desired.
Hence, $h^{-1}$ is continuous and $h$ is a homeomorphism.
\enddemo

 {\it Step} 7 (\/{\it  Completion of the proof}\/).
We return to the completion of the proof of Theorem 6.2.
Let $U_L=q_X^{-1}(W_A')\cap X$ and $U=U_L\cup U_B$. Thus, $U$ is an
open neighborhood of $Y$ in $X$. Let $U_R=U\setminus
q_X^{-1}((\xi_A')^{-1}(A\times [0,8]\times [0,+\infty)))$. 
Note that $U=U_L\cup U_R$
and 
$$U_L\cap U_R=q_X^{-1}((\xi_A')^{-1}(A\times (8,10)\times [0,+\infty))).$$
Define 
$$\varphi:U\to Y\times [0,+\infty)$$ 
as follows: 
\vglue4pt
\item{(1)} $\varphi|U_L$ is the composition 
$$U_L \mol{q_X}  W_A' \mol{\xi_A'}  A\times [0,10)\times
[0,+\infty) \moll{\pi\times\id_{[0,+\infty)}}  A\times [0,+\infty)
\mot{\subseteq} Y\times [0,+\infty).$$ 
\vglue4pt
\item{(2)} $\varphi|U_R$ is the composition 
$$U_R \mol{\varphi_B}  B\times [0,+\infty)
\moll{q_Y\times\id_{[0,+\infty)}}  Y'\times
[0,+\infty)\moll{\pi\times\id_{[0,+\infty)}} 
Y\times [0,+\infty).$$ 
\vglue4pt\noindent  
In order to verify that these definitions of $\varphi$ agree 
on the overlap, first note that
$$\multline
U_L\cap U_R = q_X^{-1}((\varphi_{AB}')^{-1}(A\times (8,10)\times
[0,+\infty)) \\
= \varphi_B^{-1}(\varphi_A^{-1}(A\times (8,10))\times [0,+\infty))
\subseteq U_{AB}.
\endmultline$$
From the definition of $q_Y$ it follows that the composition 
$$U_L\cap U_R \mol{\varphi_B}  \varphi_A^{-1}(A\times (8,10))\times
[0,+\infty) \moll{q_Y\times\id_{[0,+\infty)}}  Y'\times [0,+\infty)$$
is the composition
$$\multline
U_L\cap U_R \mol{\varphi_B}  \varphi_A^{-1}(A\times (8,10))\times
[0,+\infty) \moll{\varphi_A\times\id_{[0,+\infty)}}  
A\times (8,10)\times [0,+\infty)\\
\molm{\inclusion}  Y'\times [0,+\infty).
\endmultline$$
In turn, by the definition of $\varphi_{AB}$, this is the composition
$$U_L\cap U_R \mol{\varphi_{AB}} 
A\times (8,10)\times [0,+\infty) \molm{\inclusion}  Y'\times [0,+\infty).$$
Use the definition of $\varphi_{AB}'$ to express $\varphi_{AB}=
\varphi_{AB}'\circ q_X$.  
It follows that $\varphi\vert U_R$ on $U_L\cap U_R$ is given by 
the composition
$$
\multline
U_L\cap U_R \moll{\varphi_{AB}'\circ q_X} 
A\times (8,10)\times [0,+\infty) 
\moll{\pi\times\id_{[0,+\infty)}}  A\times [0,+\infty)\\
\molm{\inclusion}  Y\times [0,+\infty) . \endmultline$$
Using the properties of $\xi_A'$ in Step 4, this is the composition
$$\multline
U_L\cap U_R \moll{\xi_{A}'\circ q_X}   
A\times (8,10)\times [0,+\infty) 
\moll{\pi\times\id_{[0,+\infty)}}  A\times [0,+\infty)\\
\molm{\inclusion}  Y\times [0,+\infty)
\endmultline
$$
and, hence, we have the desired agreement. \pagebreak

In fact, we have shown that the compositions
$$U_L \mol{\xi_{A}'\circ q_X}  
A\times [0,10)\times [0,+\infty) 
\molm{\inclusion}  Y'\times [0,+\infty)$$
and
$$U_R \morl{(q_Y\times\id_{[0,+\infty)})\circ \varphi_B} 
Y'\times [0,+\infty)$$
agree on $U_L\cap U_R$. Hence we have a map
$$\varphi':U\to Y'\times [0,+\infty)$$
and our goal now is to show that 
$$\varphi:U \mol{\varphi'}  Y'\times [0,+\infty)
\moll{\pi\times\id_{[0,+\infty)}}  Y\times [0,+\infty)$$
has the property that its restriction
$$\varphi\vert:U\setminus Y\to Y\times (0,+\infty)$$
is an {\rm MSAF}.
This is accomplished by the following claims:

\proclaim{Claim} $\varphi'\vert:U\setminus Y\to Y'\times
(0,+\infty)$ is a {\rm MSAF}.
\endproclaim

\demo{Proof} According to the characterization in \cite{13}  
it suffices to show that the
mapping cylinder $\cyl(\varphi'\vert)$ is a homotopically
stratified space. Since this condition is a local one
it follows from the fact that $\varphi'\vert$ is locally a {\rm MSAF}.
 \enddemo

\proclaim{Claim} $\varphi\vert:U\setminus Y\to Y\times
(0,+\infty)$ is an {\rm MSAF}.
\endproclaim

\demo{Proof} First $\pi\times\id_{[0,+\infty)}$ is a {\rm MSAF} by
Claim 6.9 and [16, 4\.3]. Now combine Claim 6.12 and
 [16, 4\.5].
 \enddemo

Finally, note that $\varphi|:Y\to Y\times\{0\}$ is the
map $\pi\circ q_X=\beta$ so that Lemma 6.11 can be applied to
show that $U$ is stratum preserving homeomorphic to the teardrop
$(U\setminus Y)\cup_{\varphi|(U\setminus Y)}Y$. This shows that
$Y$ has an approximate tubular neighborhood in $X$ and completes the
proof of Theorem~6.2.
\hfill\qed

\proclaim{{C}orollary} Let $X$ be a manifold stratified space
with a compact singular set $X_{\rm sing}$ such that all nonminimal 
strata of $X$ are of dimension greater than or equal to five.
Then $X_{\rm sing}$ has an approximate tubular neighborhood in~$X$.
\endproclaim

\demo{Proof} \hskip9pt
$X_{\rm sing}$ is a singular up-set satisfying the hypothesis of\break
Theorem~6.2.
 \enddemo
\pagebreak
 
\section{Proof of the main result}
 
In this section we restate and prove the main result.

\proclaim{Theorem} Let $X$ be a manifold stratified space
with compact singular set $X_{\rm sing}$ such that all nonminimal
strata of $X$ are of dimension greater than or equal to five.
If $Y\subseteq X_{\rm sing}$ is a pure subset of $X${\rm ,}
then $Y$ has an approximate tubular neighborhood in $X$.  
\endproclaim

\demo{Proof} 
The proof is by induction on the number $n$ of strata of $X$. We
may assume that $n>0$ and that the result is true for manifold
stratified spaces with fewer than $n$ strata. Given $X$ and $Y$
as in the hypothesis, let $Z$ be the union of the strata of
$X_{\sing}\setminus Y$ and let $A=Y\cap\cl(Z)$. Then $A$ is a
pure subset of $X_{\sing}$ and $A\subseteq
{(X_{\sing})}_{\sing}$. Since $X_{\sing}$ has fewer strata than
$X$, it follows that $A$ has an approximate tubular neighborhood
in $X_{\sing}$. Say $U$ is an open neighborhood of $A$ in $X_{\sing}$ for
which there is an approximate tubular neighborhood map 
$$p:U\to A\times [0,+\infty).$$ 
By Corollary 6.14,
$X_{\sing}$ has an approximate tubular neighborhood in $X$. Say
$V$ is an open neighborhood of $X_{\sing}$ in $X$ for which there
is an approximate tubular neighborhood map 
$$q:V\to X_{\sing}\times [0,+\infty).$$ 
Note that $U\cup Y$ is open
in $X_{\sing}$. Let $U_1=U\setminus(Y\setminus A)$ so that
$U_1\cap Y=A$.
Let $W=q^{-1}((U\cup Y)\times [0,+\infty))$. Then $q|:W\to
(U\cup Y)\times [0,+\infty)$ is still a {\rm MSAF} [16, 7\.4]. 
Define 
$$\tilde p:U\cup Y\to (A\times [0,+\infty))\cup
(Y\times\{ 0\})$$ 
(where the range is a subset of $Y\times
[0,+\infty)$) by $\tilde p|U_1 =p|:U_1\to A\times
[0,+\infty)$ and $\tilde p|:Y\to Y\times\{ 0\}$ is
the identity. Note that $p|:U_1\to A\times [0,+\infty)$ is a
{\rm MSAF} because $U_1$ is a pure subset of $U$. It follows from
[16, 7\.2] that $\tilde p$ is a {\rm MSAF}.
Define 
$$r:A\times [0,+\infty)\times [0,+\infty)\to
A\times [0,+\infty)$$ 
by $r(a,s,t) = (a, s+t)$. 
It follows from [16, 4\.6] that $r$ is a stratified fibration.
Define 
$$\tilde r:[(A\times [0,+\infty))\cup Y\times\{ 0\} ]
\times [0,+\infty)\to
Y\times [0,+\infty)$$ 
by 
$$\tilde r(x,s,t) = \cases
r(x,s,t) & \text{if $ x\in A$}\\
(x,t) & \text{if $s=0$.}\\
\endcases$$
It follows from [16, 7\.2] that $\tilde r$ is a {\rm MSAF}.
Consider the composition 
$$\multline
f:W \mot{q|}  (U\cup Y)\times [0,+\infty)
\moll{\tilde p\times\id_{[0,+\infty)}}  \\
[(A\times [0,+\infty))\cup Y\times\{ 0\} ]
\times [0,+\infty) \mot{\tilde r} 
Y\times [0,+\infty).
\endmultline$$
It follows from [16, 4\.3, 4\.5]
that $f$ is a {\rm MSAF}. It is then easy to check
that $f$ is an approximate tubular neighborhood map for $Y$ in $X$.
 \enddemo

\numbereddemo{Remark}
The theorem also applies to a manifold stratified space $X$ with
noncompact singular set provided, in addition, that all the
noncompact strata are of dimension greater than or equal to five
and the one-point compactification of $X$ is a manifold
stratified space with the point at infinity constituting a new stratum.
\enddemo

\AuthorRefNames [33]
\references

 [1] 
 \name{S.\ Cappell} and \name{J.\ Shaneson},
 The mapping cone and cylinder of a stratified map, in
{\it  Prospects in Topology\/}, {\it Proc.\ Conference in Honor of
William Browder\/} (F.\ Quinn, ed.), {\it Ann. of Math.\ Studies\/} {\bf
138}, 58--66, Princeton  Univ. Press, Princeton, NJ, 1995.

[2]
 \name{T.\ A.\   Chapman},
{\it  Lectures on Hilbert Cube Manifolds\/},
 Conf. Board of the Math. Sci., 
{\it Conf.\ Series in Math\/}.\ {\bf 28}, 
 Amer.\ Math.\ Soc., 1976 
 Providence, RI, 1976.

[3] \bibline,   Approximation results in topological manifolds,
 {\it Trans.\ Amer.\ Math.\ Soc\/}.\ {\bf 262} (1980),
303--334.

[4]  \bibline, {\it Approximation Results in Hilbert Cube
Manifolds},
{\it Mem.\ Amer.\ Math.\ Soc\/}.\ {\bf 34}, No.\ 251, 1981.

[5]   \name{T.\ A.\ Chapman} and \name{S.\ Ferry},
 Approximating homotopy equivalences by homeomorphisms,
{\it Amer.\ J.\ Math\/}.\ {\bf 101} (1979), 583--607.

[6]  \name{R.\ Daverman},
{\it  Decompositions of Manifolds\/}, Academic Press, Orlando, FL, 1986.

[7]  \name{R.\ D.\ Edwards},
 TOP regular neighborhoods,
 handwritten manuscript, 1973.

[8] 
 \name{R.\ D.\ Edwards} and \name{R.\ C.\ Kirby},
Deformations of spaces of imbeddings,
{\it Ann.\ of Math\/}.\ {\bf 93} (1971), 63--88.

[9] \name{M.\ Goresky} and \name{R.\ MacPherson},
{\it  Stratified Morse Theory\/}, {\it Ergeb.\ Math.\ Grenzgeb\/}.\ {\bf
14}, Springer-Verlag, New York, 
1988.

[10] \name{B.\ Hughes},
 Approximate fibrations on topological manifolds,
{\it  Michigan Math.\ J\/}.\ {\bf 32} (1985), 167--183.

[11]     \bibline,
   Geometric topology of stratified spaces,
{\it  Electron.\ Res.\ Announc.\ Amer.\ Math.\ Soc\/}.\ 
{\bf 2} (1996),  73--81.

[12] \bibline, Stratified path spaces and fibrations,
{\it  Proc.\ Roy. Soc.\ Edinburgh Sect.\ A\/} {\bf 129} (1999), 351--384.

[13] \name{B.\ Hughes},  Stratifications of mapping cylinders,
{\it Topology Appl\/}.\ {\bf 94} (1999), 127--145.

[14] \bibline,
Stratifications of teardrops,
{\it  Fund.\ Math\/}.\  {\bf 161} (1999), 305--324.

[15]  \bibline,
Neighborhoods of strata in manifold stratified spaces,
 Vanderbilt University preprint, July 2000.

[16] \bibline,
 Products and adjunctions of manifold stratified spaces,
{\it Topology and its Appl\/}., to appear.

[17] \name{B.\ Hughes} and \name{A.\ Ranicki},
{\it Ends of Complexes\/}, {\it Cambridge Tracts in Math\/}.\
{\bf 123}, Cambridge Univ.\ Press, Cambridge, 1996.

[18] \name{B.\ Hughes, L.\ Taylor, S.\ Weinberger}, and \name{B.\ Williams},
 Neighborhoods in stratified spaces with two strata,
{\it Topology\/} {\bf 39} (2000), 873--919.

[19]  \name{B.\ Hughes, L.\ Taylor}, and \name{B.\ Williams},
 Bundle theories for topological manifolds,
{\it  Trans.\ Amer.\ Math.\ Soc\/}.\ {\bf 319} (1990), 1--65.

[20] \bibline,
Manifold approximate fibrations are approximately bundles,
{\it Forum Math\/}.\ {\bf 3} (1991), 309--325.

[21]  \name{B.\ Hughes} and \name{S.\ Weinberger},
 Surgery and stratified spaces, in
{\it  Surveys on Surgery Theory\/}, Vol.\ 2
(S.\ Cappell, A.\ Ranicki, and J.\ Rosenberg, eds.),
{\it Ann.\ of Math.\ Studies\/} {\bf 149}, 
 Princeton Univ.\ Press, 
 Princeton, NJ, 319--352.

[22] \name{J.\ Mather},
{\it  Notes on Topological Stability\/}, photocopied notes,  Harvard
Univ.,  Cambridge, MA, 1970.

[23] \bibline,
 Stratifications and mappings,
in {\it  Dynamical Systems\/}, 195--232,
{\it  Proc.\ Sympos\/}., {\it Univ.\ Bahia\/}  (Salvador,
 Brazil, 1971) (M.\ M.\ Peixoto, ed.), Academic Press, New York,
1973.

[24] \name{F.\ Quinn},
 Applications of topology with control, 598--606,
{\it Proc.\ Internat.\ Congress of Mathematicians\/} (Berkeley, CA, 1986),
 Amer.\ Math.\ Soc., Providence, RI, 1987.

[25]  \bibline,
 Homotopically stratified sets,
{\it J.\ Amer.\ Math. Soc\/}.\ {\bf 1} (1988), 441--499.

[26]  \name{C.\ Rourke} and \name{B.\ Sanderson}, 
 An embedding without a normal microbundle,
{\it  Invent.\ Math\/}.\ {\bf 3}
(1967), 293--299.

[27] \name{L.\ Siebenmann},
 Deformation  of homeomorphisms on stratified sets. I, II,  
{\it Comment.\ Math. Helv\/}.\ {\bf 47} (1972), 123--163.

[28] \name{M.\ Steinberger} and \name{J.\ West},
 Approximation by equivariant homeomorphisms, I,
{\it  Trans.\ Amer.\ Math.\ Soc\/}.\ 
{\bf 302} (1987),   297--317.

[29] \name{R.\ Thom},
 Ensembles et morphismes stratifi\'es,
{\it Bull.\ Amer.\ Math.\ Soc\/}.\ {\bf 75} (1969), 240--284. 

[30]  \name{S.\ Weinberger},
{\it  The Topological Classification of Stratified Spaces\/},
{\it  Chicago Lectures in Math\/}., Univ.\ of Chicago Press, Chicago, IL,
1994.

[31] \bibline,
 Nonlocally linear manifolds and orbifolds, 637--647,
in {\it Proc.\ Internat.\ Congress of Mathematicians\/} (Z\"urich, 
Switzerland, 1994), Birkh\"auser, Basel, 1995.

[32] \bibline,
  Microsurgery on stratified spaces, 509--521, 
in {\it Geometric Topology\/} (Athens, GA, 1993),
{\it AMS\/}/{\it IP Stud.\ Adv.\ Math\/}.\ {\bf 2.1} 
(William H.\ Kazez, ed.), Amer.\ Math.\ Soc., Providence, RI, 1997.
\endreferences

\bye

%%%%%%%%%%%%%%%%%%%%%%%%%%%%%%%%%%%%%%%%%%%%%%%%%%%%%%%%%%%%%%%%%%%%%%%

%%%%%%%%%%%%%%%%%%%%%%%%%%%%%%%%%%%%%%%%%%%%%%%%%%%%%%%%%%%%%%%%%%%%%%%%%%%%%
%%%%%%%%%%%%%%%%%%%%%%%%%%%%%%%%%%%%%%%%%%%%%%%%%%%%%%%%%%%%%%%%%%%%%%%

\Refs
\widestnumber\no{200}

\ref\no 1
\by S.\ Cappell and J.\ Shaneson
\paper The mapping cone and cylinder of a stratified map
\inbook Prospects in topology;
Proceedings of a conference in honor of William Browder
\ed F. Quinn
\bookinfo Ann. of Math. Studies\vol 138
\publ Princeton Univ. Press\publaddr Princeton\yr 1995
\pages 58--66
\endref


\ref\no 2
\by T. A. Chapman
\book Lectures on Hilbert cube manifolds
\bookinfo Conf. Board of the Math. Sci., Reg.
Conf. Series in Math. No\. 28
\publ Amer. Math. Soc. \yr 1976 
\publaddr Providence
\endref

\ref\no 3
\bysame% T. A. Chapman
\paper Approximation results in Hilbert cube manifolds
\jour Trans. Amer. Math. Soc. \vol 262 \yr 1980\pages 303--334
\endref

\ref \no 4
\bysame%   T. A. Chapman
\paper Approximation results in topological manifolds
\jour Mem. Amer. Math. Soc. 
\vol 34 \text{\rm no. 251}
\yr 1981
\endref


\ref\no 5
\by   T. A. Chapman and S. Ferry
\paper Approximating homotopy equivalences by homeomorphisms
\jour Amer. J. Math.
\vol 101
\yr 1979
\endref

\ref\no 6
\by R. Daverman
\book Decompositions of manifolds
\publ Academic Press \yr 1986
\publaddr Orlando
\endref


\ref\no 7
\by R. D. Edwards
\paper TOP regular neighborhoods
\paperinfo handwritten manuscript
\yr 1973
\endref



\ref\no 8
\by R. D. Edwards and R. C. Kirby
\paper Deformations of spaces of imbeddings
\jour Ann. of Math. (2)
\vol 93
\yr 1971
\pages 63--88
\endref

\ref\no 9
\by M. Goresky and R. MacPherson
\book Stratified Morse theory
\bookinfo Ergeb. Math. Grenzgeb (3) 14
\publ Springer-Verlag
\publaddr New York
\yr 1988
\endref



\ref\no 10 
\by B. Hughes
\paper Approximate fibrations on topological manifolds
\jour Michigan Math. J.
\vol 32
\yr 1985
\pages 167--183
\endref


\ref\no 11
\bysame%      B. Hughes
\paper   Geometric topology of stratified spaces
\jour Electron. Res. Announc. Amer. Math. Soc.
\vol 2\yr 1996\pages 73--81
\endref

\ref\no 12
\bysame
\paper Stratified path spaces and fibrations
\jour Proc. Roy. Soc. Edinburgh Sect. A
\vol 129\yr 1999\pages 351--384
\endref

\ref\no 13
\bysame
\paper Stratifications of mapping cylinders
\jour Topology Appl.
\vol 94\yr 1999\pages 127--145
\endref

\ref\no 14
\bysame 
\paper Stratifications of teardrops
\jour Fund. Math.
\vol 161 \yr 1999 \pages 305--324
\endref


\ref\no 15
\bysame %B. Hughes
\paper Neighborhoods of strata in manifold stratified spaces
\finalinfo Vanderbilt University preprint \yr July 2000
\endref

\ref\no 16
\bysame %B. Hughes
\paper Products and adjunctions of manifold stratified spaces
\jour Topology and its Appl.\finalinfo to appear
\endref

\ref \no 17
\by B. Hughes and A. Ranicki
\book Ends of complexes
\bookinfo Cambridge Tracts in Math. {\bf 123}
\publ Cambridge Univ. Press
\publaddr Cambridge 
\yr 1996
\endref


\ref\no 18
\by  B. Hughes, L. Taylor, S. Weinberger and B. Williams
\paper Neighborhoods in stratified spaces with two strata
\jour Topology
\pages 873--919\vol 39 \yr 2000
\endref


\ref\no 19
\by  B. Hughes, L. Taylor and B. Williams
\paper Bundle theories for topological manifolds
\jour Trans. Amer. Math. Soc.
\vol 319
\yr 1990
\pages 1--65
\endref


\ref\no 20
\bysame % Hughes, Taylor, Williams
\paper Manifold approximate fibrations are approximately bundles
\jour Forum Math.
\vol 3
\yr 1991
\pages 309--325
\endref

\ref \no 21
\by B. Hughes and S. Weinberger
\paper Surgery and stratified spaces
\inbook Surveys on surgery theory
\eds S. Cappell, A. Ranicki and J. Rosenberg
\vol 2\yr 2000\publ Princeton Univ. Press 
\publaddr Princeton
\pages 311--342
\endref

\ref\no 22
\by  J. Mather
\book Notes on topological stability
\publ Harvard Univ.
\publaddr Cambridge
\yr 1970
\nofrills\finalinfo (photocopied)  
\endref


\ref\no 23
\bysame % J. Mather
\paper Stratifications and mappings
\inbook Dynamical Systems
\bookinfo Proc. Sympos., Univ. Bahia, Salvador, Brazil, 1971
\ed M\. M\. Peixoto
\publ Academic Press
\publaddr New York
\yr 1973
\pages 195--232
\endref

\ref\no 24
\by F. Quinn
\paper Applications of topology with control
\inbook Proceedings of the International Congress of
Mathematicians (Berkeley, Calif., 1986)
\publ Amer. Math. Soc. 
\publaddr Providence \yr 1987 \pages 598--606
\endref


\ref\no 25
\bysame % F. Quinn
\paper Homotopically stratified sets
\jour J. Amer. Math. Soc.
\vol 1 \yr 1988 \pages 441--499
\endref

\ref\no 26
\by C. Rourke and B. Sanderson 
\paper An embedding without a normal bundle
\jour Invent. Math.
\vol 3
\yr 1967 \pages 293--299
\endref


\ref\no 27
\by L. Siebenmann
\paper Deformations of homeomorphisms on stratified sets
\jour Comment. Math. Helv.
\vol 47 \yr 1971 \pages 123--165
\endref

\ref\no 28
\by M. Steinberger and J. West
\paper Approximation by equivariant homeomorphisms. {\rom I}
\jour Trans. Amer. Math. Soc. 
\vol 302  \yr 1987 \pages 297--317
\endref

\ref\no 29
\by R. Thom
\paper Ensembles et morphismes stratifies
\jour Bull. Amer. Math. Soc.
\vol 75\yr 1969\pages 240--282
\endref

\ref\no 30
\by S. Weinberger
\book The topological classification of stratified spaces
\bookinfo Chicago Lectures in Math.
\publ Univ. Chicago Press
\publaddr Chicago
\yr 1994
\endref

\ref\no 31
\bysame% S. Weinberger
\paper   Nonlocally linear manifolds and orbifolds
\inbook Proceedings of the International Congress of
Mathematicians. Z\"urich, Switzerland 1994
\publ Birkh\"auser
\publaddr Basel
\yr 1995
\pages 637--647
\endref

\ref\no 32
\bysame% S. Weinberger
\paper   Microsurgery on stratified spaces
\inbook Geometric Topology
\ed William H. Kazez
\publ American Mathematical Society and International Press
\yr 1997
\pages 509--521
\endref

\endRefs
\enddocument